\documentclass[12pt,reqno]{amsart}

\usepackage{amssymb}

\usepackage{eucal}

\usepackage[all]{xy}

\font\sss=cmss8

\def\cK{{\mathcal K}}

\def\cO{{\mathcal O}}

\def\sD{\mbox{\sf D}}

\def\sT{\mbox{\sf T}}

\def\D{\sD}

\def\Dsmall{\mbox{\sss D}}

\def\EssIm{\operatorname{Ess.\!Im}}

\def\Hom{\operatorname{Hom}}

\def\Ker{\operatorname{Ker}}

\def\L{\operatorname{L}\!}

\def\R{\operatorname{R}\!}

\numberwithin{equation}{part}

\newtheorem{Lemma}{Lemma}[section]

\newtheorem{Theorem}[Lemma]{Theorem}
\newtheorem{Proposition}[Lemma]{Proposition}

\theoremstyle{definition}

\begin{document}

\title[Wrong way recollement]
{Wrong way recollement for schemes}

\author{Peter J\o rgensen}
\address{Department of Pure Mathematics, University of Leeds,
Leeds LS2 9JT, United Kingdom}
\email{popjoerg@maths.leeds.ac.uk}
\urladdr{http://www.maths.leeds.ac.uk/\~{ }popjoerg}






\maketitle

The notion of a recollement of triangulated categories was introduced
in \cite[sec.\ 1.4]{BBD}.  The prototypical example of a recollement
is
\[
\xymatrix
{
&&\D(Z)
    \ar[rr]^{i_*} & & 
  \D(X) 
    \ar[rr]^{j^*}
    \ar@/_1.5pc/[ll]_{i^*} \ar@/^1.5pc/[ll]^{i^!} & &
  \D(U).
    \ar@/_1.5pc/[ll]_{j_!} \ar@/^1.5pc/[ll]^{j_*} & &
}
\]
Here $X$ is a topological space equal to the union of the closed
subset $Z$ and the open complement $U$, and $\D(Z)$, $\D(X)$, and
$\D(U)$ are suitable derived categories of sheaves.  The triangulated
functors in a recollement must satisfy various conditions, most
importantly that $(i^*,i_*)$, $(i_*,i^!)$, $(j_!,j^*)$, and
$(j^*,j_*)$ are adjoint pairs.

The purpose of this note is to point out that, somewhat surprisingly,
in the case of schemes, there is also a recollement which goes the
other way.  By way of notation, if $X$ is a scheme then $\D(\cO_X)$,
the derived category of sheaves of $\cO_X$-modules, has the full
subcategory $\D(X)$ consisting of complexes with quasi-coherent
cohomology.  If $Z$ is a closed subscheme, there is also the full
subcategory $\D_Z(X)$ consisting of complexes with quasi-coherent
cohomology supported on $Z$.

\begin{Theorem}
Let $X$ be a quasi-compact quasi-separated scheme with a closed
subscheme $Z$, suppose that the open subscheme $U = X - Z$ is
quasi-compact, and write $U \stackrel{u}{\longrightarrow} X$ for the
inclusion.  Then there is a recollement
\[
\xymatrix
{
&&\D(U)
    \ar[rr]^{\R u_*} & & 
  \D(X)
    \ar[rr]
    \ar@/_1.5pc/[ll]_{\L u^*} \ar@/^1.5pc/[ll] & &
  \D_Z(X)
    \ar@/_1.5pc/[ll]_{v} \ar@/^1.5pc/[ll]^{} & &
}
\]
where $v$ is the inclusion of the full subcategory $\D_Z(X)$.
\end{Theorem}

Before proving the Theorem, let me give the following Proposition
which appears not to be stated in the literature, although it follows
easily from known results and is well known to a number of people.

\begin{Proposition}
Let $\sT$ be a compactly generated triangulated category with a set
$\cK$ which consists of compact objects and is closed under
(de)suspension.  Then there is a recollement
\[
\xymatrix
{
&&\cK^{\perp}
    \ar[rr]^{i_*} & & 
  \sT 
    \ar[rr]^{j^*}
    \ar@/_1.5pc/[ll]_{i^*} \ar@/^1.5pc/[ll]^{i^!} & &
  \langle \cK \rangle
    \ar@/_1.5pc/[ll]_{j_!} \ar@/^1.5pc/[ll]^{j_*} & &
}
\]
where $i_*$ and $j_!$ are inclusions of full subcategories.
\end{Proposition}

\noindent
Here
\[
  \cK^{\perp} = \{\, M \in \sT \,|\, 
                   \Hom_{\mbox{\sss T}}(K,M) = 0
                   \mbox{ for each $K$ in $\cK$ } \},
\]
while $\langle \cK \rangle$ is the smallest triangulated subcategory
of $\sT$ which contains $\cK$ and is closed under set indexed
coproducts; see \cite[def.\ 3.2.9]{NeemanBook}.
\medskip

{
\noindent
{\it Proof. }
It follows from the Thomason Localization Theorem, \cite[thm.\
2.1.1]{NeemanDuality}, that $\langle \cK \rangle$ is a compactly
generated triangulated category which is generated by the objects
of $\cK$.  It is clear that the inclusion functor
\[
\xymatrix
{
  \sT & &
  \langle \cK \rangle
    \ar@/_1.5pc/[ll]_{j_!}
}
\]
is a triangulated functor respecting set indexed coproducts, so by the
Brown Representability Theorem, \cite[thm.\ 4.1]{NeemanDuality}, there
is a right adjoint $j^*$ to $j_!$,
\[
\xymatrix
{
  \sT 
    \ar[rr]^{j^*} & &
  \langle \cK \rangle.
    \ar@/_1.5pc/[ll]_{j_!}
}
\]

By \cite[thm.\ 2.1.3]{NeemanDuality}, the compact objects in the
compactly generated category $\langle \cK \rangle$ are precisely the
objects of $\langle \cK \rangle$ which are compact when viewed in
$\sT$.  Hence $j_!$ sends compact objects to compact objects, and so
the right-adjoint $j^*$ respects set indexed coproducts by
\cite[thm.\ 5.1]{NeemanDuality}.  And $j^*$ is triangulated by
\cite[lem.\ 5.3.6]{NeemanBook}, so by \cite[thm.\
4.1]{NeemanDuality} again, there is a right-adjoint $j_*$ to $j^*$,
\[
\xymatrix
{
  \sT 
    \ar[rr]^{j^*} & &
  \langle \cK \rangle.
    \ar@/_1.5pc/[ll]_{j_!} \ar@/^1.5pc/[ll]^{j_*}
}
\]

As the inclusion $j_!$ is full and faithful, this is the situation
considered in \cite[prop.\ 2.7]{Miyachi} which gives a recollement
\[
\xymatrix
{
&&\Ker j^*
    \ar[rr]^{i_*} & & 
  \sT 
    \ar[rr]^{j^*}
    \ar@/_1.5pc/[ll]_{i^*} \ar@/^1.5pc/[ll]^{i^!} & &
  \langle \cK \rangle
    \ar@/_1.5pc/[ll]_{j_!} \ar@/^1.5pc/[ll]^{j_*} & &
}
\]
where $i_*$ is the inclusion of the full subcategory $\Ker j^*$.

Finally, for $N$ to be in $\Ker j^*$ means $j^* N = 0$ which holds
precisely if $\Hom_{\langle \cK \rangle}(K,j^* N) = 0$ for each $K$ in
$\cK$, because the objects of $\cK$ generate $\langle \cK \rangle$.
But
\[
  \Hom_{\langle \cK \rangle}(K,j^* N)
  \cong \Hom_{\mbox{\sss T}}(j_! K,N)
  = \Hom_{\mbox{\sss T}}(K,N),
\]
so this again holds precisely if $\Hom_{\mbox{\sss T}}(K,N) = 0$ for
each $K$ in $\cK$, that is, if $N$ is in $\cK^{\perp}$.

So I can replace $\Ker j^*$ by $\cK^{\perp}$, and this gives the
recollement of the Proposition.
\hfill $\Box$
\medskip
}

{
\noindent
{\it Proof }(of Theorem).  The category $\D_Z(X)$ is a triangulated
sub\-ca\-te\-go\-ry of $\D(X)$ which is closed under set indexed
coproducts.  Hence $\D_Z(X)$ is a triangulated category with set
indexed coproducts.  Moreover, by \cite[thm.\ 6.8]{Rouquier} there is
an object $E$ which is compact in $\D(X)$, sits in $\D_Z(X)$, and
whose (de)suspensions generate $\D_Z(X)$.  Let $\cK$ consist of all
(de)su\-spen\-si\-ons of $E$.  Then $\D_Z(X)$ is a compactly generated
triangulated category with compact generators the objects of $\cK$.

By \cite[thm.\ 2.1.2]{NeemanDuality}, it follows that $\langle \cK
\rangle_{\Dsmall_Z(X)}$, the smallest triangulated subcategory of
$\D_Z(X)$ which contains $\cK$ and is closed under set indexed
coproducts, is equal to $\D_Z(X)$ itself.  But $\langle \cK
\rangle_{\Dsmall_Z(X)}$ is clearly equal to $\langle \cK
\rangle_{\Dsmall(X)}$, the smallest triangulated subcategory of
$\D(X)$ which contains $\cK$ and is closed under set indexed
coproducts, so $\langle \cK \rangle_{\Dsmall(X)}$ is equal to
$\D_Z(X)$.

I shall simply denote $\langle \cK \rangle_{\Dsmall(X)}$ by
$\langle \cK \rangle$, so
\begin{equation}
\label{equ:a}
  \langle \cK \rangle = \D_Z(X).
\end{equation}

By \cite[prop.\ 6.7]{Spaltenstein} there is an adjoint pair of
functors
\[
\xymatrix
{
  \D(\cO_U)
    \ar[rr]_{\R u_*} & &
  \D(\cO_X).
    \ar@/_1.5pc/[ll]_{\L u^*}
}
\]
Note that since $u$ is the inclusion of an open subscheme, $u^*$ and
hence also $\L u^*$ is just restriction to $U$.  It turns out that $\L
u^*$ and $\R u_*$ restrict to the full subcategories $\D(U)$ and
$\D(X)$, and so induce an adjoint pair
\[
\xymatrix
{
  \D(U)
    \ar[rr]_{\R u_*} & &
  \D(X).
    \ar@/_1.5pc/[ll]_{\L u^*}
}
\]
This is obvious for $\L u^*$ which is just restriction to $U$.  For
$\R u_*$ it follows from \cite[thm.\ 3.3.3]{BVdB} because $u$ is a
quasi-compact quasi-separated morphism.

I will now show
\begin{equation}
\label{equ:b}
  \cK^{\perp} = \EssIm \R u_*
\end{equation}
where $\cK^{\perp}$ is taken inside $\D(X)$, while $\R u_*$ is viewed
as a functor $\D(U) \longrightarrow \D(X)$ and $\EssIm$ denotes
essential image.  To see $\supseteq$, let $M$ be in $\EssIm \R u_*$;
that is, $M \cong \R u_* N$ for some $N$ in $\D(U)$.  Then
\begin{align*}
  \Hom_{\Dsmall(X)}(K,M) 
  & \cong \Hom_{\Dsmall(X)}(K,\R u_* N) \\
  & \cong \Hom_{\Dsmall(U)}(\L u^* K,N) \\
  & \stackrel{\rm (a)}{=} 0
\end{align*}
for each $K$ in $\cK$, where (a) is because the cohomology of $K$ is
supported on $Z$ whence $\L u^* K \cong 0$.  Hence $M$ is in
$\cK^{\perp}$.

To see $\subseteq$, note first that each $M$ in $\D(X)$
determines a unit morphism
\[
  M \longrightarrow \R u_* \L u^* M.
\]
Let $I$ be a K-injective resolution of $M$ which exists by \cite[thm.\
4.5]{Spaltenstein}.  Since $\L u^*$ is just restriction, $\L u^* I
\cong u^* I$.  And $u^*$ has an exact left-adjoint $u_!$ so $u^* I$ is
also K-injective whence $\R u_* u^* I \cong u_* u^* I$.  Hence up to
isomorphism, the unit morphism is just the canonical morphism
\[
  I \longrightarrow u_* u^* I.
\]
The restriction of this to $U$ is an isomorphism because $u$ is
the inclusion of $U$ into $X$, so completing the unit morphism to a
distinguished triangle
\[
  M \longrightarrow \R u_* \L u^* M \longrightarrow C \longrightarrow,
\]
the cohomology of the cone $C$ is supported on $Z$, that is, $C$ is in
$\D_Z(X)$.

Now let $M$ be in $\cK^{\perp}$; that is,
\[
  \Hom_{\Dsmall(X)}(K,M) = 0
\]
for each $K$ in $\cK$.  Since also
\[
  \Hom_{\Dsmall(X)}(K,\R u_* \L u^* M) 
  \cong \Hom_{\Dsmall(U)}(\L u^* K,\L u^* M) = 0
\]
for each $K$ in $\cK$, the distinguished triangle implies
\[
  \Hom_{\Dsmall(X)}(K,C) = 0
\]
for each $K$ in $\cK$.  But $C$ is in $\D_Z(X)$ so this implies $C =
0$.  The distinguished triangle thus shows that the unit morphism is
an isomorphism,
\[
  M \cong \R u_* \L u^* M,
\]
so $M \cong \R u_* N$ for an $N$ in $\D(U)$.  Hence $M$ is in $\EssIm
\R u_*$ as desired.

To conclude the proof, observe that each $N$ in $\D(U)$
determines a counit morphism 
\[
  \L u^* \R u_* N \longrightarrow N.
\]
Let $J$ be a K-injective resolution of $N$.  Then $\R u_* J \cong u_*
J$ and hence $\L u^* \R u_* J \cong u^* u_* J$.  Hence up to
isomorphism, the counit morphism is just the canonical morphism
\[
  u^* u_* J \longrightarrow J
\]
which is an isomorphism because $u$ is the inclusion of an open
subscheme. 

By adjoint functor theory this implies that 
\begin{equation}
\label{equ:c}
  \mbox{$\R u_*$ is a full embedding of $\D(U)$ into $\D(X)$.}
\end{equation}

Let me now use the Proposition with $\sT = \D(X)$ and $\cK$ from
above, $\D(X)$ being compactly generated by \cite[thm.\
3.1.1(2)]{BVdB}.  The Proposition gives a recollement.  Equation
\eqref{equ:a} says that I can replace $\langle \cK \rangle$ with
$\D_Z(X)$.  Equation \eqref{equ:b} says that I can replace
$\cK^{\perp}$ with $\EssIm \R u_*$.  And equation \eqref{equ:c} says
that $\R u_*$ is a full embedding of $\D(U)$ into $\D(X)$, so I can
replace $\EssIm \R u_*$ with $\D(U)$.  This gives the recollement of
the Theorem.
\hfill $\Box$
\medskip
}

\bigskip

\noindent
{\bf Acknowledgement.}  I thank professor Leovigildo Alonso
Tarr\'{\i}o for stimulating correspondence.  The diagrams were typeset
with XY-pic.


\begin{thebibliography}{9}





\bibitem{BBD} A.\ A.\ Beilinson, J.\ Bernstein, and P.\ Deligne, {\it
Faisceaux pervers}, Ast\'{e}risque {\bf 100} (1982), 5--171 (Vol.\ 
1 of the proceedings of the conference ``Analysis and topology on
singular spaces'', Luminy, 1981).


\bibitem{BVdB}  A.\ Bondal and M.\ Van den Bergh, {\it Generators and
representability of functors in commutative and noncommutative
geometry}, Mosc.\ Math.\ J.\ {\bf 3} (2003), 1--36.
































\bibitem{Miyachi}  J.-I.\ Miyachi, {\it Localization of triangulated
categories and derived categories}, J.\ Algebra {\bf 141} (1991),
463--483.


\bibitem{NeemanDuality}  A.\ Neeman, {\it The Grothendieck duality
theorem via Bousfield's techniques and Brown representability,} J.\
Amer.\ Math.\ Soc.\ {\bf 9} (1996), 205--236.

\bibitem{NeemanBook}  \bysame, ``Triangulated categories'', Ann.\
of Math.\ Stud., Vol.\ 148, Princeton University Press, Princeton,
2001. 



\bibitem{Rouquier}  R.\ Rouquier, Dimensions of triangulated
categories, preprint (2003).  {\tt math.CT/0310134 v3}.

\bibitem{Spaltenstein}  N.\ Spaltenstein, {\it Resolutions of
unbounded complexes}, Compositio Math.\ {\bf 65} (1988), 121--154.



\end{thebibliography}
\end{document}